\documentclass[12pt,lh]{article}

\newtheorem{thm}{Theorem}[section]

\newtheorem{lem}[thm]{Lemma}
\newcounter{bean}
\newcounter{egg}

\begin{document}
\title{{\bf An Extension of the Exponential-type Error Bounds for Multiquadric and Gaussian Interpolations}}
\author{{\bf Lin-Tian Luh}\\
Department of Mathematics, Providence University\\Shalu, Taichung\\
email:ltluh@pu.edu.tw\\phone:(04)26328001 ext. 15126\\fax:(04)26324653}

\maketitle

\bigskip
{\bf Abstract.} In the 1990's exponential-type error bounds appeared in the theory of radial basis functions. For multiquadric interpolation it is $O(\lambda^{\frac{1}{d}})$ as $d\rightarrow 0$, where $\lambda$ is a constant satisfying $0<\lambda <1$. For Gaussian interpolation it is $O(C'd)^{\frac{c'}{d}}$ as $d\rightarrow 0$ where $C'$ and $c'$ are constants. In both cases the parameter $d$, called fill distance, measures the spacing of the points where interpolation occurs. This kind of error bounds is very powerful. However it only measures the difference between the approximant and approximand. Mathematicians and engineers often need to know the matching of the derivatives when dealing with partial differential equations. In this paper we extend this kind of error bounds to a form which measures the difference between the derivatives of the approximant and approximand.\\
\\
{\bf keywords}: radial basis function, error bound, multiquadric, Gaussian, derivative.

\pagenumbering{arabic}
\setcounter{section}{0}

\section{Introduction}

In this paper $h$ always denotes a continuous function on $R^{n}$ which is conditionally positive definite of order $m$. The interpolation occurs at scattered data points $(x_{j},y_{j}),\ j=1,\ldots , N$, where $X=\{ x_{1},\ldots ,x_{N}\}$ is a subset of $R^{n}$ and the $y_{j}'s$ are real or complex numbers. We are going to use as interpolant function of the form
\begin{equation}
  s(x)=p(x)+\sum_{j=1}^{N}c_{j}h(x-x_{j}),
\end{equation}
where $p(x)$ is a polynomial in $P_{m-1}$ and the $c_{j}'s$ are real or complex numbers satisfying
$$\sum_{j=1}^{N}c_{j}q(x_{j})=0$$
for all polynomials $q$ in $P_{m-1}$ and
\begin{equation}
  p(x_{i})+\sum_{j=1}^{N}c_{j}h(x_{i}-x_{j})=y_{i},\ i=1,\ldots , N.
\end{equation}
Here $P_{m-1}$ denotes the class of polynomials of $R^{n}$ of degree less than or equal to $m-1$.

A sufficient condition for a unique solution of $(1)$ is that $h$ be strictly conditionally positive definite when $X$ is a determining set for $P_{m-1}$. Recall that $X$ is said to be a determining set for $P_{m-1}$ if $p$ is in $P_{m-1}$ and $p$ vanishes on $X$ imply that $p$ is identically zero.

What's noteworthy is that in this paper both the interpolating and interpolated functions belong to a function space which was also used by Madych and Nelson in \cite{MN3} where they first introduced the exponential-type error bounds. We do not shrink or widen the function space. If $h$ is continuous and conditionally positive definite of order $m$, then it induces a native space denoted by ${\cal C}_{h,m}$ which is a semi-Hilbert space for $m>0$ and a Hilbert space for $m=0$. The definition and relevant properties of the native space can be found in \cite{MN1} and \cite{MN2}. It's not very easy to introduce this space thoroughly in limited pages. However Luh made a lucid characterization for this space in \cite{L1} and \cite{L2}. For a better understanding of this space, we suggest that the reader read \cite{L1} and \cite{L2} first.

\section{Multiquadric Interpolation}

Although Madych and Nelson made an incontrovertible contribution in the construction of the exponential-type error bounds, there are some drawbacks in their error bounds. For example, the calculation of the crucial constants in their error bounds is a big problem. Moreover, strictly speaking, their error bounds essentially apply to $R^{1}$ only. In order to make them suitable for $R^{n},\ n>1$, some modification is necessary. Recently Luh made a thorough treatment for these problems in \cite{L3} and \cite{L4}. Of course it doesn't lower the author's respect for them. However, we choose to cite Luh's error bounds in \cite{L3} and \cite{L4} since they are more complete and it will make this paper more useful.

For any set $\Omega \subseteq R^{n}$ and $X=\{ x_{1},\ldots, ,x_{N}\} \subseteq R^{n}$, let
$$d(\Omega ,X)=sup_{y\in \Omega}inf_{x\in X}|y-x| $$
be the fill distance which will be denoted by $d$ for simplicity in the following discussion. Now we need a lemma which is just Corollary2.5 of \cite{L3} with some notations simplified.
\begin{lem}
  Suppose $h$ is defined by 
$$h(x):=\Gamma(-\frac{\beta}{2})(c^{2}+|x|^{2})^{\frac{\beta}{2}},\ \beta \in R \backslash 2N_{\geq 0},\ c>0,$$
where $|x|$ is the Euclidean norm of $x$ in $R^{n}$, $\Gamma$ is the classical gamma function and $\beta ,\ c$ are real constants. Let $m=\lceil \frac{\beta}{2}\rceil $ for $\beta >0$ and $m=0$ for $\beta <0$ be its order of conditional positive definiteness. Then, given a positive number $b_{0}$, there are positive constants $d_{0}$ and $\lambda$, $0<\lambda <1$, which depend on $b_{0},\ c$ and $\beta$ for which the following is true: If $f\in {\cal C}_{h,m}$ and $s$ is the interpolant defined in (1) that interpolates $f$ on a subset $X$ of $R^{n}$, then
\begin{equation}
  |f(x)-s(x)|\leq C\cdot (\lambda)^{\frac{1}{d}}\cdot \| f\| _{h}
\end{equation}
holds for $0<d\leq d_{0}$ and all $x$ in a set $\Omega \subseteq R^{n}$, where $\Omega$ is any set in which for any point $x$, there exists a set $E$ such that $x\in E\subseteq \Omega$ where $E$ is a cube or the rotation of a cube of side $b_{0}$ satisfying the property that every subcube of $E$ of side $2d$ contains a point of $X$. Here, $\| f\| _{h}$ denotes the $h$-norm of $f$ in ${\cal C}_{h,m}$, and $C$ is a constant depending on $n,\ c$ and $\beta$.
\end{lem}
{\bf Remark.} For the calculation of $C$ and $\lambda$ in the preceding lemma we refer the reader to \cite{L3}. As for ${\cal C}_{h,m}$ and $\| f\| _{h}$, we suggest \cite{L1},\cite{L2},\cite{MN1},\cite{MN2} and \cite{MN3} for detailed introduction. However, if the reader needs a cursory understanding only, maybe \cite{MN3} is enough.
\begin{thm}
  Let $h,\ X$ and $s$ be as in Lemma2.1. For any $b_{0}>0$, let $d_{0},\ \lambda$ and $\Omega$ be as in Lemma2.1. Let $d=d(\Omega,X)$ be the fill distance, and $m$ be the order of conditional positive definiteness of $h$. Suppose $l\geq max\{ 1,m\} $ is an integer. Then
\begin{list}%
{(\alph{bean})}{\usecounter{bean} \setlength{\rightmargin}{\leftmargin}}
 \item for every $f\in {\cal C}_{h,m}$ and $|\alpha|=l$, the distribution $D^{\alpha}f$ belongs to $C(R^{n})$ and there is a constant $c_{\alpha}$ such that for all $f$ in ${\cal C}_{h,m}$,
\begin{equation}
  \| D^{\alpha}f-D^{\alpha}s\| _{\infty} \leq c_{\alpha}\| f\| _{h}\ ;
\end{equation}
 \item for any $y\in \Omega $ such that the closed ball $\overline{B}(y,\delta )\subseteq \Omega $ for some $\delta>0$, and any $f\in {\cal C}_{h,m}$, if $0<|\alpha |<l$, then 
\begin{equation}
  |D^{\alpha }f(y)-D^{\alpha }s(y)|\leq  C_{|\alpha |}M_{0}^{1-\frac{|\alpha |}{l}}\cdot \overline{M}_{l}^{\frac{|\alpha |}{l}}
\end{equation}
where $C_{|\alpha |}$ is a constant independent of $f$, $M_{0}=C\| f\| _{h}\lambda ^{\frac{1}{d}}$ as in (3), and $\overline{M}_{l}=max(M_{l},\ M_{0}l!\delta ^{-l})$ with $M_{l}:=C'\| f\| _{h}$ where $C'=l!\sum _{|\alpha |=l}\frac{c_{\alpha}}{\alpha !}$.
\end{list}

\end{thm}
Proof. To prove (a), note that $s\in {\cal C}_{h,m}$ \cite{MN1} and hence $f-s\in {\cal C}_{h,m}$. Then (a) follows immediately from Proposition4.1. of \cite{MN2} and the discussion preceding it. What's difficult is to show (b).

For $y\in \Omega $ with $\overline{B}(y,\delta )\subseteq \Omega $ for some $\delta >0$, $t\in [-\delta ,\delta ]$ and $u\in R^{n}$ with $|u|=1$, let $\psi (t)=f(y+tu)-s(y+tu)$. Then 
$$\psi^{(k)}(t)=k!\sum _{|\alpha |=k}\frac{u^{\alpha}}{\alpha !}D^{\alpha}(f-s)(y+tu)$$
for $0<k\leq l$. By Proposition4.1 of \cite{MN2}, $\| \psi^{(l)}\| _{\infty}\leq C'\| f-s\| _{h}$ with $C'=l!\sum _{|\alpha|=l}c_{\alpha}/\alpha!$. The discussion of p.223 of \cite{MN2} tells us that $\| f-s\| _{h}\leq \| f\| _{h}$. Thus $\| \psi^{(l)}\| _{\infty }\leq C'\| f\| _{h}$. From (3) we also have a bound for $|\psi(t)|$ on the interval $[-\delta,\delta]$. By the result of Gorny in p.139 of \cite{Mi} we get for $0<k<l$,
$$|\psi^{(k)}(0)|\leq 16(2e)^{k}M_{0}^{1-\frac{k}{l}}\cdot \overline{M}_{l}^{\frac{k}{l}}$$
where $M_{0}=C\| f\| _{h}\lambda^{\frac{1}{d}}$ and $\overline{M}=max(M_{l},\ M_{0}l!\delta^{-l})$ with $M_{l}=C'\| f\| _{h}$ where $C'=l!\sum _{|\alpha|=l}c_{\alpha}/\alpha!$. This implies
$$\left| \sum _{|\alpha|=k}\frac{u^{\alpha}}{\alpha!}D^{\alpha}(f-s)(y)\right| \leq \frac{1}{k!}16(2e)^{k}M_{0}^{1-\frac{k}{l}}\overline{M}_{l}^{\frac{k}{l}}.$$
Since $u$ is any unit vector in $R^{n}$, we have
$${sup}_{|u|=1}\left| \sum_{|\alpha|=k}\frac{u^{\alpha}}{\alpha!}D^{\alpha}(f-s)(y)\right| \leq \frac{1}{k!}16(2e)^{k}M_{0}^{1-\frac{k}{l}}\overline{M}_{l}^{\frac{k}{l}}$$
for every $y\in \Omega$ with $\overline{B}(y,\delta)\subseteq \Omega$ for some $\delta >0$. Now, let $V_{k}$ be the space of vectors $v=(v_{\alpha})_{|\alpha|=k}$. Then
$$|v|_{k}=sup_{|u|=1}\left| \sum_{|\alpha|=k}\frac{u^{\alpha}}{\alpha!}v_{\alpha}\right| $$
is a norm for $V_{k}$. Any two norms on a finite-dimensional space are equivalent. Therefore this norm is equivalent to the Euclidean norm. This gives that for each $|\alpha|=k$,
$$\left| D^{\alpha}(f-s)(y)\right| \leq \frac{\rho_{k}}{k!}16(2e)^{k}M_{0}^{1-\frac{k}{l}}\overline{M}_{l}^{\frac{k}{l}}$$
for some constant $\rho_{k}>0$. Let $C_{k}=\frac{\rho_{k}}{k!}16(2e)^{k}$. Then (b) follows.\hspace{2cm}\ \ $\sharp$\\
\\
{\bf Remark:} Note that when $d$ is small, $\overline{M}_{l}=M_{l}$ and (5) becomes
\begin{eqnarray*}
  |D^{\alpha}f(y)-D^{\alpha}s(y)| & \leq & C_{|\alpha|}\cdot \{ C\| f\| _{h}\lambda^{\frac{1}{d}}\} ^{1-\frac{|\alpha|}{l}}\cdot \{ C'\| f\| _{h}\} ^{\frac{|\alpha|}{l}}\\
                                  & = & C_{|\alpha|}\cdot C^{1-\frac{|\alpha|}{l}}\cdot {C'}^{\frac{|\alpha|}{l}}\cdot \| f\| _{h}\cdot \lambda^{\frac{1}{d}(1-\frac{|\alpha|}{l})}.
\end{eqnarray*}
Thus we have
\begin{equation}
  |D^{\alpha}f(y)-D^{\alpha}s(y)| \leq C_{|\alpha|}\cdot C^{1-\frac{|\alpha|}{l}}\cdot {C'}^{\frac{|\alpha|}{l}}\cdot \| f\| _{h}\cdot \lambda^{\frac{1}{d}(1-\frac{|\alpha|}{l})}
\end{equation}
for $0<|\alpha|<l$ whenever $\overline{B}(y,\delta)\subseteq \Omega$ for some $\delta>0$ and $d$ is small.

If $d$ is large, $\overline{M}_{l}=M_{0}l!\delta^{-l}$ and (5) becomes
\begin{eqnarray*}
  |D^{\alpha}f(y)-D^{\alpha}s(y)| & \leq & C_{|\alpha|}\cdot M_{0}^{1-\frac{|\alpha|}{l}}\cdot \{ M_{0}l!\delta^{-l}\}                                             ^{\frac{|\alpha|}{l}}\\
                                  & = & C_{|\alpha|}\cdot M_{0}\cdot (l!)^{\frac{|\alpha|}{l}}\cdot \delta^{-|\alpha|}\\
                                  & = & C_{|\alpha|}\cdot C\cdot \| f\| _{h}\cdot \lambda^{\frac{1}{d}}\cdot (l!)^{\frac{|\alpha|}{l}}\cdot \delta^{-|\alpha|}. 
\end{eqnarray*}
Thus we have
\begin{equation}
  |D^{\alpha}f(y)-D^{\alpha}s(y)| \leq C_{|\alpha|}\cdot C\cdot (l!)^{\frac{|\alpha|}{l}}\cdot \delta^{-|\alpha|}\cdot \| f\| _{h}\cdot \lambda^{\frac{1}{d}}
\end{equation}
for $0<|\alpha|<l$ whenever $\overline{B}(y,\delta)\subseteq \Omega$ for some $\delta>0$ and $d$ is large.

From (6) and (7) we find that the effect of differentiation appears essentially only when $d$ is small.

\section{Gaussian Interpolation}

\begin{lem}
  Let $h(x)=e^{-\beta|x|^{2}},\ \beta>0$, be the Gaussian function in $R^{n}$. Then given a positive number $b_{0}$, there are positive constants $d_{0},\ g$ and $G$ which depend on $b_{0}$ and $\beta$ for which the following is true: If $f\in {\cal C}_{h,m}$ and $s$ is the interpolant defined in (1) that interpolates $f$ on a subset $X$ of $R^{n}$, then
\begin{equation}
  |f(x)-s(x)| \leq \Delta''(Gd)^{\frac{g}{d}}\| f\| _{h}
\end{equation}
holds for $0<d\leq d_{0}$ and all x in a set $\Omega \subseteq R^{n}$, where $\Omega$ is any set in which any point x is contained in a subset $E$ of $\Omega$, the set $E$ being a cube or the rotation of a cube of side $b_{0}$ satisfying the property that every subcube of $E$ of side $2d$ contains a point of $X$. Here, $\| f\| _{h}$ denotes the $h$-norm of $f$ in ${\cal C}_{h,m}$, and $\Delta''$ is a constant depending on n.
\end{lem}
Proof. This is just Corollary1 of \cite{L4}.\\
\\
{\bf Remark}: The calculation of $\Delta'',\ G$ and $g$ is quite complicated. The interested readers are referred to \cite{L4} for a detailed discussion.
\begin{thm}
  Let $h,\ X$ and $s$ be as in Lemma3.1. For any $b_{0}>0$, let $d_{0},\ g$ and $G$ be as in Lemma3.1 as well. Suppose $\Omega \subseteq R^{n}$ and for any $x\in \Omega$ there is a set $E$ such that $x\in E\subseteq \Omega$ where $E$ is a cube or the rotation of a cube of side $d_{0}$ satisfying the property that every subcube of $E$ of side $2d,\ 0<d\leq d_{0}$, contains a point of $X$. Let $d=d(\Omega,X)$ be the fill distance and $l\geq 1$ be any integer. Then
\begin{list}%
{(\alph{egg})}{\usecounter{egg}
   \setlength{\rightmargin}{\leftmargin}}
\item for every $f\in {\cal C}_{h,m}$ and $|\alpha|=l$, the distribution $D^{\alpha}f$ belongs to $C(R^{n})$ and there is a       constant $c_{\alpha}$ such that for all $f$ in ${\cal C}_{h,m}$, 
\begin{equation}
  \| D^{\alpha}f-D^{\alpha}s\| _{\infty}\leq c_{\alpha}\| f\| _{h}\ ;
\end{equation} 
\item for any $y\in \Omega$ such that the closed ball $\overline{B}(y,\delta)\subseteq \Omega$ for some $\delta >0$, and any $f\in {\cal C}_{h,m}$ if $\ 0<|\alpha| <l$, then 
\begin{equation}
  |D^{\alpha}f(y)-D^{\alpha}s(y)| \leq C_{|\alpha|}M_{0}^{1-\frac{|\alpha|}{l}}\cdot \overline{M}_{l}^{\frac{|\alpha|}{l}}
\end{equation} 
where $C_{|\alpha|}$ is a constant independent of $f$, $M_{0}=\Delta''(Gd)^{\frac{g}{d}}\| f\| _{h}$ as in (8), and $\overline{M}_{l}=max(M_{l},\ M_{0}l!\delta^{-l})$ with $M_{l}=C'\| f\| _{h}$ where $C'=l!\sum _{|\alpha|=l}\frac{C_{\alpha}}{\alpha !}$.
\end{list}
\end{thm}
Proof. (a) is similar to Theorem2.2 with the only difference that the Gaussian function is conditionally positive definite of order $m=0$ and $c_{\alpha}$ is defined in another way because it corresponds to a different radial function $h$.

The proof of (b) is also similar to Theorem2.2 with a different definition of $M_{0}$. We omit it too.\hspace{9cm}\ \ \  $\sharp$\\
\\
{\bf Remark}: Note that when $d$ is small, $\overline{M}_{l}=M_{l}$ and (10) becomes
\begin{eqnarray*}
  |D^{\alpha}f(y)-D^{\alpha}s(y)| & \leq & C_{|\alpha|}\{ \Delta'' (Gd)^{\frac{g}{d}}\| f\| _{h}\} ^{1-\frac{|\alpha|}{l}}\cdot \{ C'\| f\| _{h}\} ^{\frac{|\alpha|}{l}}\\
                                  & = & C_{|\alpha|}\cdot {\Delta''}^{1-\frac{|\alpha|}{l}}\cdot {C'}^{\frac{|\alpha|}{l}}\cdot \| f\| _{h}\cdot (Gd)^{\frac{g}{d}(1-\frac{|\alpha|}{l})}.
\end{eqnarray*}
Thus we have
\begin{equation}
  |D^{\alpha}f(y)-D^{\alpha}s(y)| \leq C_{|\alpha|}\cdot {\Delta''}^{1-\frac{|\alpha|}{l}}\cdot {C'}^{\frac{|\alpha|}{l}}\cdot \| f\| _{h}\cdot (Gd)^{\frac{g}{d}(1-\frac{|\alpha|}{l})}
\end{equation}
for $0<|\alpha|<l$ whenever $\overline{B}(y,\delta)\subseteq \Omega$ for some $\delta>0$ and $d$ is small.

If $d$ is large, $\overline{M}=M_{0}l!\delta^{-l}$ and (10) becomes
\begin{eqnarray*}
  |D^{\alpha}f(y)-D^{\alpha}s(y)| & \leq & C_{|\alpha|}\cdot M_{0}^{1-\frac{|\alpha|}{l}}\cdot \{ M_{0}l!\delta^{-l}\} ^{\frac{|\alpha|}{l}}\\
                                  & = & C_{|\alpha|}\cdot M_{0}\cdot (l!)^{\frac{|\alpha|}{l}}\cdot \delta^{-|\alpha|}\\
                                  & = & C_{|\alpha|}\cdot \Delta'' (Gd)^{\frac{g}{d}}\cdot \| f\| _{h}\cdot (l!)^{\frac{|\alpha|}{l}}\cdot \delta^{-|\alpha|}. 
\end{eqnarray*}
Thus we have
\begin{equation}
  |D^{\alpha}f(y)-D^{\alpha}s(y)| \leq C_{|\alpha|}\cdot \Delta''\cdot (l!)^{\frac{|\alpha|}{l}}\cdot \delta^{-|\alpha|}\cdot \| f\| _{h}\cdot (Gd)^{\frac{g}{d}}
\end{equation}
for $0<|\alpha| <l$ whenever $\overline{B}(y,\delta)\subseteq\Omega$ for some $\delta>0$ and $d$ is large.

From (11) and (12) we find that, as in the case of multiquadrics, the effect of differentiation appears essentially only when $d$ is small.

Finally, we'd like to point out the main goal we are pursuing. In recent years mathematicians noticed the importance of the matching of smoothness. In order to obtain a good approximation, the derivatives should also be approximated. This is why we present this paper.\\

\end{document}